\documentclass[twocolumn]{IEEEtran}

\usepackage{graphics,epsfig,amsmath,amsfonts,amssymb,cite,times,colortbl,subfigure,floatflt,wrapfig,enumerate,tikz,url,xcolor}
\usepackage{algorithmic,eqparbox,comment}
\usepackage{bm}
\usepackage{booktabs}

\newcommand{\thickhline}{\noalign{\hrule height 1.0pt}}
\newtheorem{definition}{Definition}
\newtheorem{example}{Example}

\newcommand{\kpr}[1]{\textsuperscript{\textcircled{#1}}}    
\newcommand{\ten}[1]{\bmcal{#1}} 
\DeclareMathAlphabet\bmcal{OMS}{cmsy}{b}{n}

\newcommand{\argmin}[1]{\underset{#1}{\operatorname{argmin}}}
\begin{document}

\title{Tensor Computation: A New Framework for High-Dimensional Problems in EDA}
\author{Zheng Zhang, Kim~Batselier, Haotian~Liu, Luca Daniel and Ngai~Wong
\thanks{Z. Zhang and L. Daniel are with Department of Electrical Engineering and Computer Science, Massachusetts Institute of Technology (MIT), Cambridge, MA. E-mails: \{z\_zhang, luca\}@mit.edu}
\thanks{K. Batselier and N. Wong are with Department of Electrical and Electronic Engineering, the University of Hong Kong. E-mails: \{kimb, nwong\}@eee.hku.hk}
\thanks{H. Liu is with Cadence Design Systems, Inc. San Jose, CA. E-mail: haotian@cadence.com}
 }
 \markboth{Accepted by IEEE TRANSACTIONS ON COMPUTER-AIDED DESIGN OF INTEGRATED CIRCUITS AND SYSTEMS, ~Vol. ~XX, No.~XX,~XX~2016}{ZHANG \MakeLowercase{\textit{et al.}}: Tensor Computation}

\IEEEspecialpapernotice{(Invited Keynote Paper)}

\maketitle
\begin{abstract}
Many critical EDA problems suffer from the curse of dimensionality, i.e. the very fast-scaling computational burden produced by large number of parameters and/or unknown variables. This phenomenon may be caused by multiple spatial or temporal factors (e.g. 3-D field solvers discretizations and multi-rate circuit simulation), nonlinearity of devices and circuits, large number of design or optimization parameters (e.g. full-chip routing/placement and circuit sizing), or extensive process variations (e.g. variability/reliability analysis and design for manufacturability). The computational challenges generated by such high dimensional problems are generally hard to handle efficiently with traditional EDA core algorithms that are based on matrix and vector computation. This paper presents ``tensor computation'' as an alternative general framework for the development of efficient EDA algorithms and tools. A tensor is a high-dimensional generalization of a matrix and a vector, and is a natural choice for both storing and solving efficiently high-dimensional EDA problems. This paper gives a basic tutorial on tensors, demonstrates some recent examples of EDA applications (e.g., nonlinear circuit modeling and high-dimensional uncertainty quantification), and suggests further open EDA problems where the use of tensor computation could be of advantage.

\end{abstract}

\section{Introduction}
\label{sec:intro}

\subsection{Success of Matrix \& Vector Computation in EDA Hystory}
The advancement of fabrication technology and the development of Electronic Design Automation (EDA) are two engines that have been driving the progress of semiconductor industries. The first integrated circuit (IC) was invented in 1959 by Jack Kilby. However, until the early 1970s designers could only handle a small number of transistors manually. The idea of EDA, namely designing electronic circuits and systems automatically using computers, was proposed in the 1960s. Nonetheless, this idea was regarded as science fiction until SPICE~\cite{SPICE1} was released by UC Berkeley. Due to the success of SPICE, numerous EDA algorithms and tools were further developed to accelerate various design tasks, and designers could design large-scale complex chips without spending months or years on labor-intensive work. 

The EDA area indeed encompasses a very large variety of diverse topics, e.g., hardware description languages, logic synthesis, formal verification. This paper mainly concerns computational problems in EDA.  Specifically, we focus on modeling, simulation and optimization problems, whose performance heavily relies on effective numerical implementation. Very often, numerical modeling or simulation core tools are called repeatedly by many higher-level EDA tools such as design optimization and system-level verification. Many efficient matrix-based and vector-based algorithms have been developed to address the computational challenges in EDA. Here we briefly summarize a small number of examples among the numerous research results. 

In the context of circuit simulation, modified nodal analysis~\cite{MNA1} was proposed to describe the dynamic network of a general electronic circuit. Standard numerical integration and linear/nonlinear equation solvers (e.g., Gaussian elimination, LU factorization, Newton's iteration) were implemented in the early version of SPICE~\cite{SPICE1}. Driven by communication IC design, specialized RF simulators were developed for periodic steady-state~\cite{Kundert:90, Aprille:ieeeProc,Aprille:TCAS, Kundert:89, Petzold:81} and noise~\cite{demir2000phase} simulation. Iterative solvers and their parallel variants were further implemented to speed up large-scale linear~\cite{kozhaya2002multigrid,Chen:2001,feng2008multigrid} and nonlinear circuit simulation~\cite{Ricardo:1995,liu2015gpu}. In order to handle process variations, both Monte Carlo~\cite{MCintro, SingheeR10} and stochastic spectral methods~\cite{zzhang_cicc2014, zzhang_iccad2013,zzhang:huq_tcad, twweng:optsEx,zzhang:tcad2013, zzhang:tcas2_2013, Pulch:2011_1, Wang:2004,Vrudhula:2006, Rufuie2014, manfredi:tcas2014}) were investigated to accelerate stochastic circuit simulation.

Efficient models were developed at almost every design level of hierarchy. At the process level, many statistical and learning algorithms were proposed to characterize manufacturing process variations~\cite{variation2008, yu2014remembrance,zhang2011virtual}. At the device level, a huge number of physics-based (e.g., BSIM~\cite{chauhan2012bsim} for MOSFET and RLC interconnect models) and math-based modeling frameworks were reported and implemented. Math-based approaches are also applicable to circuit and system-level problems due to their generic formulation. They start from a detailed mathematical description [e.g., a partial differential equation (PDE) or integral equation describing device physics~\cite{FastCap, FastHenrry, pFFT} or a dynamic system describing electronic circuits] or some measurement data, then generate compact models by model order reduction~\cite{OCP:98,  PDS:03, JR:99,LP:05, P:03, G:11,ZLWFW:12,BondD09,RW:03} or system identification~\cite{Gustavsen1999, grivet2004passivity, coelho2004convex, bond2010sysid}. These techniques were further extended to problems with design parameters or process uncertainties~\cite{DOLLW:04,DOLLW:02,sou2008quasi,BD:05, BondD07,MD:10,ferranti2010guaranteed, villena2010spare, silveira2006resampling}.

Thanks to the progress of numerical optimization~\cite{boyd2004convex,vandenberghe1996semidefinite,bertsekas1999nonlinear}, a lot of algorithmic solutions were developed to solve EDA problems such as VLSI placement~\cite{Shahookar1991}, routing~\cite{cong1996performance}, logic synthesis~\cite{de1994synthesis} and analog/RF circuit optimization~\cite{Gielen:optimization1990,cai2011optimization}. Based on design heuristics or numerical approximation, the performance of many EDA optimization engines were improved. For instance, in analog/RF circuit optimization, posynomial or polynomial performance models were extracted to significantly reduce the number of circuit simulations~\cite{hershenson2001optimal, li2004robust, xu2005opera}. 

\subsection{Algorithmic Challenges and Motivation Examples}
Despite the success in many EDA applications, conventional matrix-based and vector-based algorithms have certain intrinsic limitations when applied to problems with high dimensionality. These problems generally involve an extremely large number of unknown variables or require many simulation/measurement samples to characterize a quantity of interest. Below we summarize some representative motivation examples among numerous EDA problems:

\begin{itemize}
\item \textbf{Parameterized 3-D Field Solvers.} Lots of devices are described by PDEs or integral equations~\cite{FastCap, FastHenrry, pFFT} with $d$ spatial dimensions. With $n$ discretization elements along each spatial dimension (i.e., $x$-, $y$- or $z$-direction), the number of unknown elements is approximately $N$=$n^d$ in a finite-difference or finite-element scheme. When $n$ is large (e.g. more than thousands and often millions), even a fast iterative matrix solver with $O(N)$ complexity cannot handle a 3-D device simulation. If design parameters (e.g. material properties) are considered and the PDE is further discretized in the parameter space, the computational cost quickly extends beyond the capability of existing matrix- or vector-based algorithms.

\item \textbf{Multi-Rate Circuit Simulation.} Widely separated time scales appear in many electronic circuits (e.g. switched capacitor filters and mixers), and they are difficult to simulate using standard transient simulators. Multi-time PDE solvers~\cite{Jaijeet:PDE2001} reduce the computational cost by discretizing the differential equation along $d$ temporal axes describing different time scales. Similar to a 3-D device simulator, this treatment may also be affected by the curse of dimensionality. Frequency-domain approaches such as multi-tone harmonic balance~\cite{Melville:mt1995, Carvalho:mt1998} may be more efficient for some RF circuits with $d$ sinusoidal inputs, but their complexity also becomes prohibitively high as $d$ increases.

\item  \textbf{Probabilistic Noise Simulation.} When simulating a circuit influenced by noise, some probabilistic approaches (such as those based on Fokker-Planck equations~\cite{bonnin2013phase}) compute the joint density function of its $d$ state variables along the time axis. In practice, the $d$-variable joint density function must be finely discretized in the $d$-dimensional space, leading to a huge computational cost. 

\item \textbf{Nonlinear or Parameterized Model Order Reduction.} The curse of dimensionality is a long-standing challenge in model order reduction. In multi-parameter model order reduction~\cite{DOLLW:04,DOLLW:02,villena2010spare, silveira2006resampling}, a huge number of moments must be matched, leading to a huge-size reduced-order model. In nonlinear model order reduction based on Taylor expansions or Volterra series~\cite{JR:99,LP:05,P:03}, the complexity is an exponential function of the highest degree of Taylor or Volterra series. Therefore, existing matrix-based algorithms can only capture low-order nonlinearity.

\item  \textbf{Design Space Exploration.} Consider a classical design space exploration problem: optimize the circuit performance (e.g., small-signal gain of an operational amplifier) by choosing the best values of $d$ design parameters (e.g. the sizes of all transistors). When the performance metric is a strongly nonlinear and discontinuous function of design parameters, sweeping the whole parameter space is possibly the only feasible solution. Even if a small number of samples are used for each parameter, a huge number of simulations are required to explore the whole parameter space.

\item \textbf{Variability-Aware Design Automation.} Process variation is a critical issue in nano-scale chip design. Capturing the complex stochastic behavior caused by process uncertainties can be a data-intensive task. For instance, a huge number of measurement data points are required to characterize accurately the variability of device parameters~\cite{variation2008, yu2014remembrance, zhang2011virtual}. In circuit modeling and simulation, the classical stochastic collocation algorithm~\cite{Pulch:2011_1,Wang:2004,Vrudhula:2006} requires many simulation samples in order to construct a surrogate model. Although some algorithms such as compressed sensing~\cite{xli2011, zhang2011virtual} can reduce measurement or computational cost, lots of hidden data information cannot be fully exploited by matrix-based algorithms.

\end{itemize}

\subsection{Toward Tensor Computations?}
In this paper we argue that one effective way to address the above challenges is to utilize tensor computation. Tensors are high-dimensional generalizations of vectors and matrices. Tensors were developed well over a century ago, but have been mainly applied in physics, chemometrics and psychometrics~\cite{tensorreview}. Due to their high efficiency and convenience in representing and handling huge data arrays, tensors are only recently beginning to be successfully applied in many engineering fields, including (but not limited to) signal processing~\cite{Cichocki2015}, big data~\cite{vervliet2014}, machine learning and scientific computing. Nonetheless, tensors still seem a relatively unexplored and unexploited concept in the EDA field.

The goals and organization of this paper include:
\begin{itemize}
	\item Providing a hands-on ``primer'' introduction to tensors and their basic computation techniques (Section~\ref{sec:tensor} and appendices), as well as the most practically useful techinques such as tensor decomposition (Section \ref{sec:tensordecomp}) and tensor completion (Section \ref{sec:tensorcompl});
	\item Summarizing, as guiding examples, a few recent tensor-based EDA algorithms, including progress in high-dimensional uncertainty quantification (Section~\ref{sec:applications}) and nonlinear circuit modeling and simulation (Section~\ref{sec:nonlinear});
	\item Suggesting some theoretical and application open challenges in tensor-based EDA (Sections~\ref{sec:future_app} and \ref{sec:future}) in order to stimulate further research contributions.
\end{itemize}

\section{Tensor Basics}
\label{sec:tensor}
This section reviews some basic tensor notions and operations necessary for understanding the key ideas in the paper. Different fields have been using different conventions for tensors. Our exposition will try to use one of the most popular and consistent notations. 

\subsection{Notations and Preliminaries}
\begin{figure}[t]
\centering 
\includegraphics[width=2.0in]{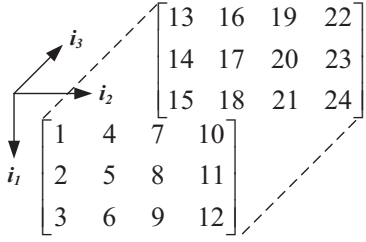}
\caption{An example tensor $\ten{A}\in\mathbb{R}^{3\times 4\times 2}$.}
\label{fig:tensorexample}
\end{figure}

We use boldface capital calligraphic letters (e.g. $\ten{A}$) to denote tensors, boldface capital letters (e.g. $\bm{A}$) to denote matrices, boldface letters (e.g. $\bm{a} $) to denote vectors, and roman (e.g. $a$) or Greek (e.g. $\alpha$) letters to denote scalars.

\textbf{Tensor.} A tensor is a high-dimensional generalization of a matrix or vector. A vector $\bm{a}\in \mathbb{R}^{n}$ is a 1-way data array, and its $i$th element $a_i$ is specified by the index $i$. A matrix $\bm{A} \in \mathbb{R}^{n_1\times n_2}$ is a 2-way data array, and each element $a_{i_1 i_2}$ is specified by a row index $i_1$ and a column index $i_2$. By extending this idea to the high-dimensional case $d\geq 3$, a tensor $\ten{A}\in\mathbb{R}^{n_1\times n_2\times\cdots\times n_d}$ represents a $d$-way data array, and its element $a_{i_1 i_2 \cdots i_d}$ is specified by $d$ indices. Here, the positive integer $d$ is also called the order of a tensor. Fig.~\ref{fig:tensorexample} illustrates an example $3\times 4\times 2$ tensor.

\subsection{Basic Tensor Arithmetic}
\label{subsec:tensorarithmetic}
\definition{\textbf{Tensor inner product.}}
The inner product between two tensors $\ten{A},\ten{B} \in \mathbb{R}^{n_1 \times \cdots \times n_d}$ is defined as
$$
\langle \ten{A},\ten{B} \rangle \;=\; \sum_{i_1,i_2,\ldots,i_d} a_{i_1\cdots i_d} b_{i_1 \cdots i_d}.
$$
As norm of a tensor $\ten{A}$, it is typically convenient to use the Frobenius norm $||\ten{A}||_F := \sqrt{\langle \ten{A},\ten{A} \rangle}$.

\definition{\textbf{Tensor $\bm{k}$-mode product.}}
The $k$-mode product $\ten{B}=\ten{A} {\times_k}\, \bm{U}$ of a tensor $\ten{A}\in\mathbb{R}^{n_1\times\cdots \times n_k\times\cdots\times n_d}$ with a matrix $\bm{U}\in\mathbb{R}^{p_k\times n_k}$ is defined by%
\begin{align}
b_{i_1\cdots i_{k-1} j i_{k+1} \cdots i_d}=\sum\limits_{i_k=1}^{n_k}  u_{j i_k} a_{i_1\cdots i_k\cdots i_d},%
\label{eqn:kmode}
\end{align}
and $\ten{B}\in\mathbb{R}^{n_1\times\cdots \times n_{k-1}\times p_k\times n_{k+1}\times\cdots\times n_d}$.

\definition{\textbf{$\bm{k}$-mode product shorthand notation.}}
The multiplication of a $d$-way tensor $\ten{A}$ with the matrices $\bm{U}^{(1)},\ldots,\bm{U}^{(d)}$ along each of its $d$ modes respectively is
\begin{align*}
[\![\ten{A};\bm{U}^{(1)}, \ldots, \bm{U}^{(d)}]\!] \triangleq \ten{A} \times_1 \bm{U}^{(1)} \times_2 \cdots \times_d \bm{U}^{(d)}.
\end{align*}
When $\ten{A}$ is diagonal with all 1's on its diagonal and 0's elsewhere, then $\ten{A}$ is omitted from the notation, e.g. $[\![\bm{U}^{(1)}, \ldots, \bm{U}^{(d)}]\!]$.

\definition{\textbf{Rank-1 tensor.}} A rank-1 $d$-way tensor can be written as the outer product of $d$ vectors
\begin{align}
\ten{A} = \bm{u}^{(1)} \circ \bm{u}^{(2)} \circ \cdots \circ \bm{u}^{(d)} = [\![\bm{u}^{(1)}, \ldots, \bm{u}^{(d)}]\!],
\label{eqn:outerprod}
\end{align}
where $\bm{u}^{(1)} \in \mathbb{R}^{n_1},\ldots,\bm{u}^{(d)}\in \mathbb{R}^{n_d}$. The entries of $\ten{A}$ are completely determined by $a_{i_1 i_2\cdots i_d}=u^{(1)}_{i_1}u^{(2)}_{i_2}\cdots u^{(d)}_{i_d}$.

Some additional notations and operations are introduced in Appendix~\ref{appd:notations}. The applications in Sections~\ref{sec:applications} and \ref{sec:nonlinear} will make it clear that the main problems in tensor-based EDA applications are either computing a tensor decomposition or solving a tensor completion problem. Both of them will now be discussed in order.

\section{Tensor Decomposition}
\label{sec:tensordecomp}

\subsection{Computational Advantage of Tensor Decompositions.}
The number of elements in a $d$-way tensor is $n_1n_2\cdots n_d$, which grows very fast as $d$ increases. Tensor decompositions compress and represent a high-dimensional tensor by a smaller number of factors. As a result, it is possible to solve high-dimensional problems (c.f. Sections V to VII) with a lower storage and computational cost.  Table~\ref{table:tdecomp} summarizes the storage cost of three mainstream tensor decompositions in order to intuitively show their advantage. State-of-the-art implementations of these methods can be found in \cite{TTB_Software,tensorlab3,tttoolbox}. 
Specific examples are for instance: 
\begin{itemize}
\item
While the hidden layers of a neural network could consume almost all of the memory in a server, using a canonical or tensor-train decomposition instead results in an extraordinary compression(\cite{novikov2015,lebedev2014}) by up to a factor of $200,000$.
\item
High-order models describing nonlinear dynamic systems can also be significantly compressed using tensor decomposition as will be shown in details in Section~\ref{sec:nonlinear}.
\item
High-dimensional integration and convolution are long-standing challenges in many engineering fields (e.g. computational finance and image processing). These two problems can be written as the inner product of two tensors, and while a direct computation would have a complexity of $O(n^d)$, using a low-rank canonical or tensor-train decomposition, results in an extraordinarily lower $O(nd)$ complexity~\cite{rakhuba2015fast}. 
\end{itemize}
In this section we will briefly discuss the most popular and useful tensor decompositions, highlighting advantages of each.

\subsection{Canonical Polyadic Decomposition}
\label{subsec:cpd}
\begin{table}[t]
\centering
\caption{Storage costs of mainstream tensor decomposition approaches.}
\label{table:tdecomp} 
\begin{tabular}{@{}lrr@{}}
Decomposition & Elements to store & Comments \\ \midrule
\\
Canonical Polyadic~\cite{candecomp,harshman1970fpp} & {$ndr$}  & see Fig.~\ref{fig:tensorouterprod} \\
\\
Tucker~\cite{tuckerreview} & {$r^d+ndr$}  & see Fig.~\ref{fig:tuckerdecomp} \\
\\
 Tensor Train~\cite{ivanTT} & {$n(d-2)r^2+2nr$}  & see Fig.~\ref{fig:tt} 
\end{tabular}
\end{table}

\begin{figure}[t]
\centering
\includegraphics[width=3.3in]{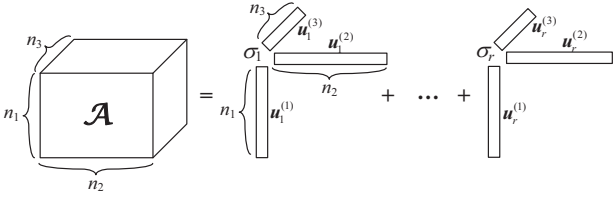}
\caption{Decomposing $\ten{A}$ into the sum of $r$ rank-$1$ outer products.}
\label{fig:tensorouterprod}
\end{figure}

\textbf{Polyadic Decomposition}. A polyadic decomposition expresses a $d$-way tensor as the sum of $r$ rank-$1$ terms:
\begin{align}
\ten{A}=\sum\limits_{i=1}^{r} \sigma_i \, \bm{u}_i^{(1)} \circ \cdots \circ \bm{u}_i^{(d)}
 = [\![\ten{D};\bm{U}^{(1)}, \cdots, \bm{U}^{(d)}]\!].
\label{eqn:outersum}%
\end{align}%
The subscript $i$ of the unit-norm $\bm{u}_i^{(1)}$ vectors indicates a summation index and not the vector entries. The $\bm{u}_i^{(k)}$ vectors are called the mode-$k$ vectors. Collecting all vectors of the same mode $k$ in matrix $\bm{U}^{(k)}\in \mathbb{R}^{n_k \times r}$, this decomposition is rewritten as the $k$-mode products of matrices $\{ \bm{U}^{(k)}\}_{k=1}^d$ with a cubical diagonal tensor $\ten{D} \in \mathbb{R}^{r \times \cdots  \times r}$ containing all the $\sigma_i$ values. Note that we can always absorb each of the scalars $\sigma_i$ into one of the mode vectors, then write $\ten{A}=[\![\bm{U}^{(1)},\cdots,\bm{U}^{(d)}]\!]$. 

\example The polyadic decomposition of a 3-way tensor is shown in Fig.~\ref{fig:tensorouterprod}.

\textbf{Tensor Rank}. The minimum $r:=R$ for the equality (\ref{eqn:outersum}) to hold is called the \textit{tensor rank} which, unlike the matrix case, is in general NP-hard to compute~\cite{Hastad1990}. 

\textbf{Canonical Polyadic Decomposition (CPD)}. The corresponding decomposition with the minimal $R$ is called the \textit{canonical polyadic decomposition} (CPD). It is also called \textit{Canonical Decomposition} (CANDECOMP)~\cite{candecomp} or \textit{Parallel Factor} (PARAFAC)~\cite{harshman1970fpp} in the literature. A CPD is unique, up to scaling and permutation of the mode vectors, under mild conditions. A classical uniqueness result for 3-way tensors is described by Kruskal \cite{Kruskal1977}. These uniqueness conditions do not apply to the matrix case\footnote{Indeed, for a given matrix decomposition $\bm{A}=\bm{U}\bm{V}$ and any nonsingular matrix $\bm{T}$ we have that $\bm{A}=\bm{U}\bm{T}\bm{T^{-1}}\bm{V}$. Only by adding sufficient conditions (e.g. orthogonal or triangular factors) the matrix decomposition can be made unique. Remarkably, the CPD for higher order tensors does not need any such conditions to ensure its uniqueness.}.

The computation of a polyadic decomposition, together with two variants are discussed in Appendix~\ref{app:cpvariant}.

\subsection{Tucker Decomposition}
\begin{figure}[t]
\centering
\includegraphics[width=3.2in]{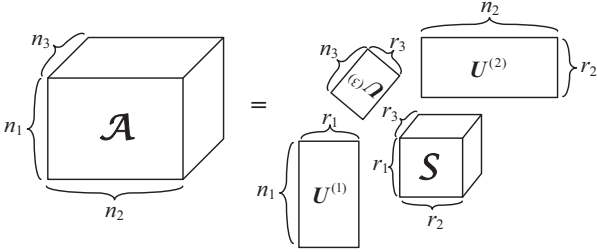}
\caption{The Tucker decomposition decomposes a 3-way tensor $\ten{A}$ into a core tensor $\ten{S}$ and factor matrices $\bm{U}^{(1)},\bm{U}^{(2)},\bm{U}^{(3)}$.}
\label{fig:tuckerdecomp}
\end{figure}

\textbf{Tucker Decomposition}. Removing the constraint that $\ten{D}$ is cubical and diagonal in \eqref{eqn:outersum} results in
\begin{align}
   \label{eqn:tucker}
   \ten{A} &=  \ten{S} \times_{1} \bm{U}^{(1)} \times_{2} \bm{U}^{(2)}\cdots \times_{d} \bm{U}^{(d)}\\
& =[\![\ten{S};\bm{U}^{(1)}, \bm{U}^{(2)},\ldots, \bm{U}^{(d)}]\!] \nonumber
\end{align}
with the factor matrices $\bm{U}^{(k)} \in \mathbb{R}^{n_k \times r_k}$ and a core tensor $\ten{S} \in \mathbb{R}^{r_1 \times r_2 \times \cdots \times r_d}$. The Tucker decomposition can significantly reduce the storage cost when $r_k$ is (much) smaller than $n_k$. This decomposition is illustrated in Fig.~\ref{fig:tuckerdecomp}. 

\textbf{Multilinear Rank}. The minimal size $(r_1,r_2,\ldots,r_d)$ of the core tensor $\ten{S}$ for \eqref{eqn:tucker} to hold is called the multilinear rank of $\ten{A}$, and it can be computed as $r_1=\textrm{rank}(\ten{A}_{(1)}),\ldots,r_d=\textrm{rank}(\ten{A}_{(d)})$. Note that $\ten{A}_{(k)}$ is a matrix obtained by reshaping (see Appendix~\ref{appd:notations}) $\ten{A}$ along its $k$th mode. For the matrix case we have that $r_1=r_2$, i.e., the row rank equals the column rank. This is not true anymore when $d \geq 3$. 

\textbf{Tucker vs. CPD}. The Tucker decomposition can be considered as an expansion in rank-1 terms that is not necessarily canonical, while the CPD does not necessarily have a minimal core. This indicates the different usages of these two decompositions: the CPD is typically used to decompose data into interpretable mode vectors while the Tucker decomposition is most often used to compress data into a tensor of smaller size. Unlike the CPD, the Tucker decomposition is in general not unique\footnote{One can always right-multiply the factor matrices $\bm{U}^{(k)}$ with any nonsingular matrix $\bm{T}^{(k)}$ and multiply the core tensor $\ten{S}$ with their inverses ${\bm{T}^{(k)}}^{-1}$.
This means that the subspaces that are defined by the factor matrices $\bm{U}^{(i)}$ are invariant while the bases in these subspaces can be chosen arbitrarily.}.

A variant of the Tucker decomposition, called high-order singular value decomposition (SVD) or HOSVD, is summarized in Appendix~\ref{app:hosvd}.

\subsection{Tensor Train Decomposition}

\textbf{Tensor Train (TT) Decomposition}. A tensor train decomposition~\cite{ivanTT} represents a $d$-way tensor $\ten{A}$ by two $2$-way tensors and $(d-2)$ $3$-way tensors. Specifically, each entry of $\ten{A} \in \mathbb{R}^{n_1\times \cdots \times n_d}$ is expressed as
\begin{equation}
a_{i_1 i_2 \cdots i_d} = \ten{G}^{(1)}_{i_1} \, \ten{G}^{(2)}_{i_2} \cdots \ten{G}^{(d)}_{i_d},
\label{eq:TTdecomp}
\end{equation}
where $\ten{G}^{(k)} $$\in $$\mathbb{R}^{r_{k-1} \times n_k \times r_k} $ is the $k$-th core tensor, $r_0=r_d=1$, and thus $\ten{G}^{(1)}$ and $\ten{G}^{(d)}$ are matrices. The vector $(r_0, r_1, \cdots, r_d)$ is called the \textit{tensor train rank}. Each element of the core $\ten{G}^{(k)}$, denoted as $g^{(k)}_{\alpha_{k-1} i_k \alpha_{k+1}}$ has three indices. By fixing the $2$nd index $i_k$, we obtain a matrix $\ten{G}^{(k)}_{i_k}$ (or vector for $k=1$ or $k=d$).

\begin{figure}[t]
\centering
\includegraphics[width=3.2in]{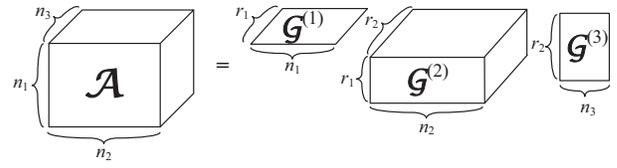}
\caption{The Tensor Train decomposition decomposes a 3-way tensor $\ten{A}$ into two matrices $\ten{G}^{(1)},\ten{G}^{(3)}$ and a $3$-way tensor $\ten{G}^{(2)}$.}
\label{fig:tt}
\end{figure}

\textbf{Computing Tensor Train Decompositions}. Computing a tensor train decomposition consists of doing $d-1$ consecutive reshapings and low-rank matrix decompositions. An advantage of tensor train decomposition is that a quasi-optimal approximation can be obtained with a given error bound and with an automatic rank determination~\cite{ivanTT}.

\subsection{Choice of Tensor Decomposition Methods}
Canonical and tensor train decompositions are preferred for high-order tensors since the their resulting tensor factors have a low storage cost linearly dependent on $n$ and $d$. For some cases (e.g., functional approximation), a tensor train decomposition is preferred due to a unique feature, i.e., it can be implemented with cross approximation~\cite{ttcross} and {\it without} knowing the whole tensor. This is very attractive, because in many cases obtaining a tensor element can be expensive. Tucker decompositions are mostly applied to lower-order tensors due to their storage cost of $O(r^d)$, and are very useful for finding the dominant subspace of some modes such as in data mining applications.

\section{Tensor Completion (or Recovery)}
\label{sec:tensorcompl}

Tensor decomposition is a powerful tool to reduce storage and computational cost, however most approaches need a whole tensor {\it a-priori}. In practice, obtaining each element of a tensor may require an expensive computer simulation or non-trivial hardware measurement. Therefore, it is necessary to estimate a whole tensor based on only a small number of available elements. This can be done by tensor completion or tensor recovery. This idea finds applications in many fields. For instance in biomedical imaging, one wants to reconstruct the whole magnetic resonance imaging data set based on a few measurements. In design space exploration, one may only have a small number of tensor elements obtained from circuit simulations, while all other sweeping samples in the parameter space must be estimated.

\subsection{Ill-Posed Tensor Completion/Recovery}  
Let ${\cal I}$ include all indices for the elements of $\ten{A}$, and its subset $\Omega$ holds the indices of some available tensor elements. A projection operator $\mathbb{P}_{\Omega}$ is defined for $\ten{A}$:
\begin{equation}
\label{tensor_project}
\ten{B}=\mathbb{P}_{\Omega}\left({\ten{A}}\right) \; \Leftrightarrow\;  b_{i_1\cdots i_d} = \left\{ \begin{array}{l}
 a_{i_1\cdots i_d} ,\;{\rm{if}}\; {i_1\cdots i_d} \in  {\Omega} \\
 0 ,\;{\rm{otherwise}}.
 \end{array} \nonumber \right.
 \end{equation}
In tensor completion, one wants to find a tensor $\ten{X}$ such that it matches $\ten{A}$ for the elements specified by $\Omega$:
 \begin{equation}
\label{tensor_point_match}
\| \mathbb{P}_{\Omega}\left(\ten{X} -{\ten{A}}\right)\|_F^2 =0.
\end{equation}
This problem is \textbf{ill-posed}, because any value can be assigned to $x_{i_1\cdots i_d}$ if $i_1\cdots i_d \notin \Omega$.

\subsection{Regularized Tensor Completion} 
Regularization makes the tensor completion problem well-posed by adding constraints to (\ref{tensor_point_match}). Several existing ideas are summarized below.
\begin{itemize}
\item \textbf{Nuclear-Norm Minimization.} This idea searches for the minimal-rank tensor by solving the problem:
\begin{equation}
\label{eq:tc_convex}
\min_{\ten{X} } { \| \ten{X}\|_{*} } \;\;\;\;\; {\rm s. } {\rm t.} \; \mathbb{P}_{\Omega}\left(\ten{X}\right)= \mathbb{P}_{\Omega} \left(\ten{A}\right) .
\end{equation}
The nuclear norm of a matrix is the sum of all singular values, but the nuclear norm of a tensor does not have a rigorous or unified definition. In~\cite{Gandy:2011, Liu:2013}, the tensor nuclear norm $\| \ten{X}\|_{*} $ is heuristically approximated using the weighted sum of matrix nuclear norms of $\ten{X}_{(k)}$'s for all modes. This heuristic makes~(\ref{eq:tc_convex}) convex, and its optimal solution can be computed by available algorithms~\cite{Douglas:1956, Gabay:1976}. Note that in~(\ref{eq:tc_convex}) one has to compute a full tensor $\ten{X}$, leading to an exponential complexity with respect to the order $d$.

\item \textbf{Approximation with Fixed Ranks.} Some techniques compute a tensor $\ten{X}$ by fixing its tensor rank. For instance, one can solve the following problem
 \begin{align}
\label{tensor_rankr}
& \min_{\ten{X}} \; {\| \mathbb{P}_{\Omega}\left(\ten{X} -{\ten{A}}\right)\|_F^2 }  \nonumber \\
 & {\rm s.}\; {\rm t.} \;\;\;{\rm multilinear}\; {\rm rank} (\ten{X}) =(r_1, \ldots, r_d)
\end{align}
with $\ten{X}$ parameterized by a proper low-multilinear rank factorization. Kresner et al.~\cite{Kressner:2014} computes the higher-order SVD representation using Riemannian optimization~\cite{Absil:2008}. In~\cite{Holtz:2012}, the unknown $\ten{X}$ is parameterized by some tensor-train factors. The low-rank factorization significantly reduces the number of unknown variables. However, how to choose an optimal tensor rank still remains an open question.

\item \textbf{Probabilistic Tensor Completion.} In order to automatically determine the tensor rank, some probabilistic approaches based on Bayesian statistics have been developed. Specifically, one may treat the tensor factors as unknown random variables assigned with proper prior probability density functions to enforce low-rank properties. This idea has been applied successfully to obtain polyadic decomposition~\cite{Rai:2014, Zhao:2015} and Tucker decomposition~\cite{Zhao:arxiv2015} from incomplete data with automatic rank determination.

\item \textbf{Low-Rank and Sparse Constraints.} In some cases, a low-rank tensor $\ten{A}$ may have a sparse property after a linear transformation. Let $\bm{z}=[z_1, \ldots, z_m]$ with $z_k=  \left\langle \ten{A}, \ten{W}_k \right\rangle$, one may find that many elements of $\bm{z}$ are close to zero. To exploit the low-rank and sparse properties simultaneously, the following optimization problem~\cite{zzhang:spi2016,zzhang:cpmt2016} may be solved:
 \begin{align}
\label{tensor_lrsp}
& \min_{\ten{X}} \; {\frac{1}{2} \| \mathbb{P}_{\Omega}\left(\ten{X} -{\ten{A}}\right)\|_F^2 +\lambda \sum\limits_{k=1}^m |\left\langle \ten{X}, \ten{W}_k\right\rangle | }  \nonumber \\
& {\rm s.} \; {\rm t.} \; \;\;{\rm multilinear}\; {\rm rank} (\ten{X}) =(r_1, \ldots, r_d).
\end{align}
In signal processing, $\bm{z}$ may represent the coefficients of multidimensional Fourier or wavelet transforms. In uncertainty quantification, $\bm{z}$ collects the coefficients of a generalized polynomial-chaos expansion. The formulation (\ref{tensor_lrsp}) is generally non-convex, and locating its global minimum is non-trivial.
\end{itemize}

\subsection{Choice of Tensor Recovery Methods}
Low-rank constraints have proven to be a good choice for instance in signal and image processing (e.g., MRI reconstruction)~\cite{Gandy:2011,Liu:2013,Kressner:2014,Zhao:2015}. Both low-rank and sparse properties may be considered for high-dimensional functional approximation (e.g., polynomial-chaos expansions)~\cite{zzhang:spi2016}. Nuclear-norm minimization and probabilistic tensor completion are very attractive in the sense that tensor ranks can be automatically determined, however they are not so efficient or reliable for high-order tensor problems. It is expensive to evaluate the nuclear norm of a high-order tensor. Regarding probabilistic tensor completion, implementation experience shows that many samples may be required to obtain an accurate result.

\section{Applications in Uncertainty Quantification}
\label{sec:applications}
Tensor techniques can advance the research of many EDA topics due to the ubiquitous existence of high-dimensional problems in the EDA community, especially when considering process variations. This section summarizes some recent progress on tensor-based research in solving high-dimensional uncertainty quantification problems, and could be used as guiding reference for the effective employment of tensors in other EDA problems.

\subsection{Uncertainty Quantification (UQ)}
\label{subsec:UQ}
Process variation is one of the main sources causing yield degradation and chip failures. In order to improve chip yield, efficient stochastic algorithms are desired in order to simulate nano-scale designs. The design problems are generally described by complex differential equations, and they have to be solved repeatedly in traditional Monte-Carlo simulators.

Stochastic spectral methods have emerged as a promising candidate due to their high efficiency in EDA applications~\cite{zzhang_cicc2014, zzhang_iccad2013,zzhang:huq_tcad, twweng:optsEx,zzhang:tcad2013, zzhang:tcas2_2013, Pulch:2011_1, Wang:2004,Vrudhula:2006, Rufuie2014, manfredi:tcas2014}. Let the random vector $\boldsymbol{\xi} \in \mathbb{R}^d$ describe process variation. Under some assumptions, an output of interest (e.g., chip frequency) $y(\boldsymbol{\xi})$ can be approximated by a truncated generalized polynomial-chaos expansion~\cite{gPC2002}:
\begin{equation}
\label{surrogate_i}
y\left (\boldsymbol{\xi}\right )\approx \sum\limits_{|\boldsymbol{\alpha}|=0}^p {c_{\boldsymbol{\alpha}} \Psi _{\boldsymbol{\alpha}}  (\boldsymbol{\xi})}.
\end{equation} \normalsize
 Here $\{\Psi _{\boldsymbol{\alpha}}(\boldsymbol{\xi})\}$ are orthonormal polynomial basis functions; the index vector $\boldsymbol{\alpha} \in \mathbb{N}^d$ indicates the polynomial order, and its element-wise sum $|\boldsymbol{\alpha}|$ is bounded by $p$. The coefficient $c_{\boldsymbol{\alpha}}$ can be computed by
\begin{equation}
\label{c_project}
c_{\boldsymbol{\alpha}}= \mathbb{E} \left( \Psi_{\boldsymbol{\alpha}} ( \boldsymbol{\xi} ) y ( \boldsymbol{\xi} )\right) 
\end{equation}
where $\mathbb{E}$ denotes expectation.

\textbf{Main Challenge}. Stochastic spectral methods become inefficient when there are many random parameters, because evaluating $c_{\boldsymbol{\alpha}}$ involves a challenging $d$-dimensional numerical integration. In high-dimensional cases, Monte Carlo was regarded more efficient than stochastic spectral methods. However, we will show that with tensor computation, stochastic spectral methods can outperform Monte Carlo for some challenging UQ problems.

\begin{figure}[t]
	\centering
		\includegraphics[width=2.0in]{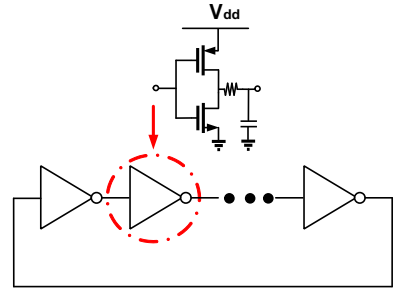} 
\caption{Schematic of a multistage CMOS ring oscillator (with 7 inverters). }
	\label{fig:ring}
\end{figure}

\begin{figure*}[t]
	\centering
		\includegraphics[width=5.2in]{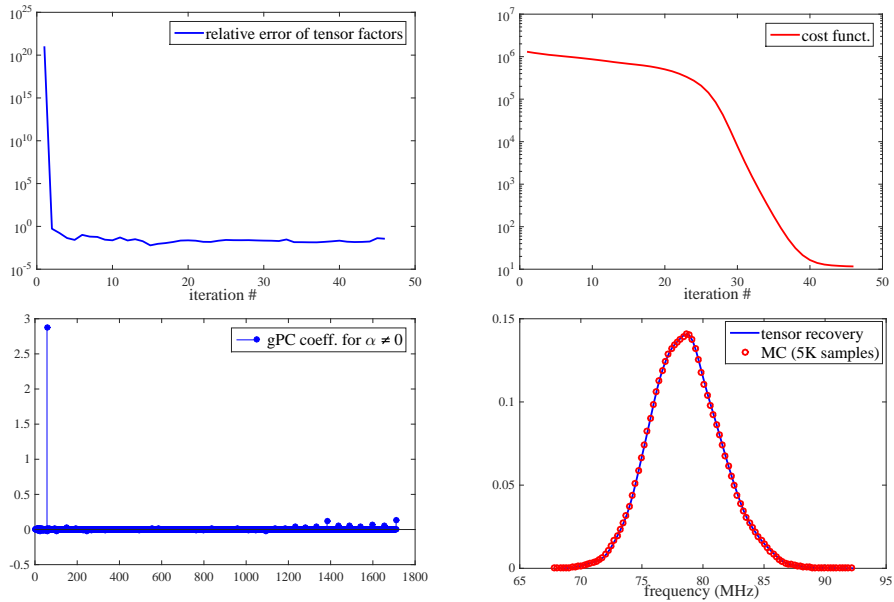} 
\caption{Tensor-recovery results of the ring oscillator. Top left: relative error of the tensor factors; top right: decrease of the cost function in (\ref{tensor_lrsp}); bottom left: sparsity of the obtained polynomial-chaos expansion; bottom right: obtained density function v.s. Monte Carlo using $5000$ samples.}
	\label{fig:ring_results}
\end{figure*}

\subsection{High-D Stochastic Collocation by Tensor Recovery} 
\textbf{Problem Description.} In stochastic collocation~\cite{col:2005, Ivo:2007, Nobile:2008}, (\ref{c_project}) is evaluated by a quadrature rule. For instance, with $n_j$ integration samples and weights~\cite{Golub:1969} properly chosen for each element of $\boldsymbol{\xi}$, $c_{\boldsymbol{\alpha}}$ can be evaluated by
\begin{equation}
\label{c_tensor}
c_{\boldsymbol{\alpha}}=\left \langle  \ten{Y}, \ten{W}_{\boldsymbol{\alpha}} \right \rangle.
\end{equation}
Here both $\ten{Y}$ and $\ten{W}_{\boldsymbol{\alpha}}$ are tensors of size $n_1\times \cdots \times n_d$. The rank-1 tensor $\ten{W}_{\boldsymbol{\alpha}}$ only depends on $\Psi_{\boldsymbol{\alpha}} ( \boldsymbol{\xi} )$ and quadrature weights, and thus is easy to compute. Obtaining $\ten{Y}$ exactly is almost impossible because it has the values of $y$ at all integration samples. Instead of computing all elements of $\ten{Y}$ by $n_1n_2 \cdots n_d$ numerical simulations, we estimate $\ten{Y}$ using only a small number of (say, several hundreds) simulations. As shown in compressive sensing~\cite{xli2011}, the approximation (\ref{surrogate_i}) usually has sparse structures, and thus the low-rank and sparse tensor completion model (\ref{tensor_lrsp}) can be used. Using tensor recovery, stochastic collocation may require only a few hundred simulations, thus can be very efficient for some high-dimensional problems.

\textbf{Example~\cite{zzhang:spi2016,zzhang:cpmt2016}.} The CMOS ring oscillator in Fig.~\ref{fig:ring} has $57$ random parameters describing threshold voltages, gate-oxide thickness, and effective gate length/width. Since our focus is to handle high dimensionality, all parameters are assumed mutually independent. We aim to obtain a $2$nd-order polynomial-chaos expansion for its frequency by repeated periodic steady-state simulations. Three integration points are chosen for each parameter, leading to $3^{57}\approx 1.6 \times 10^{27}$ samples to simulate in standard stochastic collocation. Advanced integration rules such as sparse grid~\cite{sparsegrids} still needs over $6000$ simulations. As shown in Table~\ref{table:ring_cost}, with tensor completion (\ref{tensor_lrsp}), the tensor representing $3^{57}$ solution samples can be well approximated by using only $500$ samples. As shown in Fig.~\ref{fig:ring_results}, the optimization solver converges after $46$ iterations, and the tensor factors are computed with less than $1\%$ relative errors; the obtained model is very sparse, and the obtained density function of the oscillator frequency is very accurate. 
\begin{table} [t]
\caption{Comparison of simulation cost for the ring oscillator, using different kinds of stochastic collocation.} 
\centering 
\begin{tabular}{l r r r} 
method & tensor product & sparse grid & tensor completion \\  \thickhline 
total samples & $ 1.6\times 10^{27}$ & $6844$ & $500$ \\ 
\end{tabular} 
\label{table:ring_cost}
\end{table} 

\textbf{Why Not Use Tensor Decomposition?} Since $\ten{Y}$ is not given {\it a priori}, neither CPD nor Tucker decomposition is feasible here. For the above example our experiments show that tensor train decomposition requires about $10^5$ simulations to obtain the low-rank factors with acceptable accuracy, and its cost is even higher than Monte Carlo.

\subsection{High-D Hierarchical UQ with Tensor Train} 
\begin{figure}[t]
	\centering
		\includegraphics[width=3.0in]{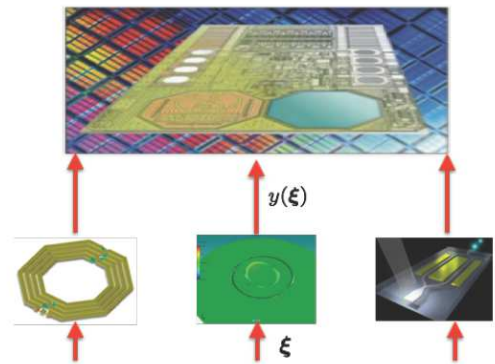} 
\caption{Hierarchical uncertainty quantification. The stochastic outputs of bottom-level components/devices are used as new random inputs for upper-level uncertainty analysis. }
	\label{fig:huq}
\end{figure}

\textbf{Hierarchical UQ.} In a hierarchical UQ framework, one estimates the high-level uncertainty of a large system that consists of several components or subsystems by applying stochastic spectral methods at different levels of the design hierarchy. Assume that several polynomial-chaos expansions are given in the form (\ref{surrogate_i}), and each $y\left (\boldsymbol{\xi}\right )$ describes the output of a component or subsystem. In Fig.~\ref{fig:huq}, $y\left (\boldsymbol{\xi}\right )$ is used as a new random input such that the system-level simulation can be accelerated by ignoring the bottom-level variations $\boldsymbol{\xi}$. However, the quadrature samples and basis functions of $y$ are unknown, and one must compute such information using a 3-term recurrence relation~\cite{Golub:1969}. This requires evaluating the following numerical integration with high accuracy:
\begin{equation}
\label{int_huq}
\mathbb{E}\left( g\left(y\left (\boldsymbol{\xi}\right ) \right )\right)= \left \langle \ten{G}, \ten {W} \right \rangle,
\end{equation}
where the elements of tensors $\ten{G}$ and $ \ten{W} \in \mathbb{R}^{\hat {n}_1 \times \cdots \times \hat{n}_d}$ are $g(y(\boldsymbol{\xi}_{i_1\cdots i_d} ))$ and $w_1^{i_1}\cdots w_d^{i_d}$, respectively. Note that $\boldsymbol{\xi}_{i_1\cdots i_d}$ and $w_1^{i_1}\cdots w_d^{i_d}$ are the $d$-dimensional numerical quadrature samples and weights, respectively.

\textbf{Choice of Tensor Decompositions.} We aim to obtain a low-rank representation of $\ten{Y}$, such that $\ten{G}$ and $\mathbb{E}\left( g\left(y\left (\boldsymbol{\xi}\right ) \right )\right)$ can be computed easily. Due to the extremely high accuracy requirement in the 3-term recurrence relation~\cite{Golub:1969}, tensor completion methods are not feasible. Neither canonical tensor decomposition nor Tucker decomposition is applicable here, as they need the whole high-way tensor $\ten{Y} $ before factorization. Tensor-train decomposition is a good choice, since it can compute a high-accuracy low-rank representation without knowing the whole tensor $\ten{Y}$; therefore, it was used in \cite{zzhang:huq_tcad} to accelerate the 3-term recurrence relation and the subsequent hierarchical UQ flow.

\begin{figure*}[t]
	\centering 		\includegraphics[width=5.5in]{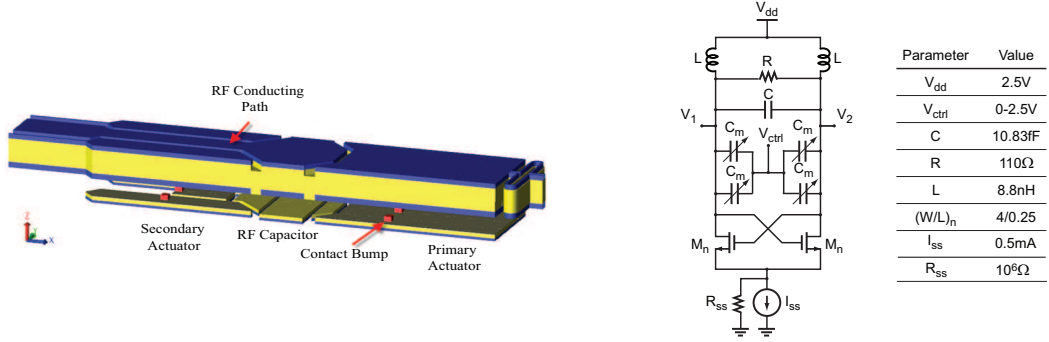} 
\caption{Left: the schematic of a MEMS switch acting as capacitor $C_{\rm m}$, which has 46 process variations; right: an oscillator using 4 MEMS switches as capacitors (with $184$ random parameters in total). }
	\label{fig:memsOsc}
\end{figure*}

\begin{figure}[t]
	\centering
		\includegraphics[width=3.1in]{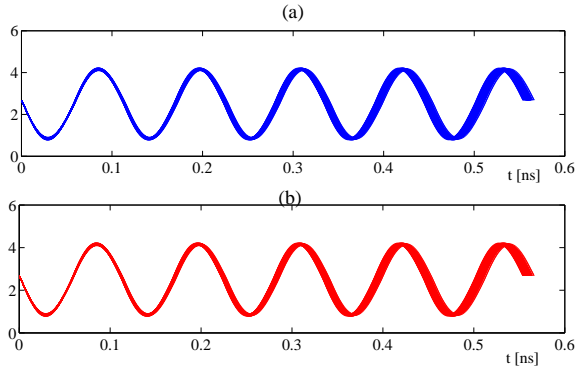} 
\caption{Realization of the steady-state waveforms for the oscillator in Fig.~\ref{fig:memsOsc}. Top: tensor-based hierarchical approach; bottom: Monte Carlo.}
	\label{fig:wave_mc}
\end{figure}

\textbf{Example.} The tensor-train-based flow has been applied to the oscillator with four MEMS capacitors and $184$ random parameters shown in Fig.~\ref{fig:memsOsc}, which previously could only be solved using random sampling approaches. In~\cite{zzhang:huq_tcad}, a sparse generalized polynomial-chaos expansion was first computed as a stochastic model for the MEMS capacitor $y(\boldsymbol{\xi})$. The discretization of $y(\boldsymbol{\xi})$ on a $46$-dimensional integration grid was represented by a tensor $\ten{Y}$ (with $9$ integration points along each dimension), then $\ten{Y}$ was approximated by a tensor train decomposition. After this approximation, (\ref{int_huq}) was easily computed to obtain the new orthonormal polynomials and quadrature points for $y$. Finally, a stochastic oscillator simulator~\cite{zzhang:tcas2_2013} was called at the system level using the newly obtained basis functions and quadrature points. As shown in Table~\ref{table:huq}, this circuit was simulated by the tensor-train-based hierarchical approach in only $10$ min in MATLAB, whereas Monte Carlo with 5000 samples required more than 15 hours~\cite{zzhang:huq_tcad}. The variations of the steady-state waveforms from both methods are almost the same, cf. Fig.~\ref{fig:wave_mc}.
\begin{table} [t]
\caption{Simulation time of the MEMS-IC co-design in Fig.~\ref{fig:memsOsc}} 
\centering 
\begin{tabular}{l r  r} 
method & Monte Carlo & proposed~\cite{zzhang:huq_tcad}\\  \thickhline 
total samples & $ 15.4$ hours& $10$ minutes\\
\end{tabular} 
\label{table:huq}
\end{table} 

\section{Applications in Nonlinear Circuit Modeling}
\label{sec:nonlinear}
Nonlinear devices or circuits must be well modeled in order to enable efficient system-level simulation and optimization. Capturing the (possibly high) nonlinearity can result in high-dimensional problems. Fortunately, the multiway nature of a tensor allows the easy capturing of high-order nonlinearities of analog, mixed-signal circuits and in MEMS design. 

\subsection{Nonlinear Modeling and Model Order Reduction}
\label{subsec:MS}
Similar to the Taylor expansion, it is shown in ~\cite{WLWW:12,LP:05, P:03,HW:13,ZW:15} that many nonlinear dynamical systems can be approximated by expanding the nonlinear terms around an equilibrium point, leading to the following ordinary differential equation
\begin{align}
\bm{\dot x}=\bm{Ax}+\bm{B} \bm{x}\kpr{2}+\bm{C}\bm{x}\kpr{3} + \bm{D}(\bm{u} \otimes \bm{x}) +\bm{Eu},
\label{eq:ode}
\end{align}
where the state vector $\bm{x}(t)\in\mathbb{R}^n$ contains the voltages and/or currents inside a circuit network, and the vector $\bm{u}(t)\in\mathbb{R}^m$ denotes time-varying input signals. The $\bm{x}\kpr{2},\bm{x}\kpr{3}$ notation refers to repeated Kronecker products (cf.~Appendix \ref{appd:notations}). The matrix $\bm{A}\in\mathbb{R}^{n\times n}$ describes linear behavior, while the matrices $\bm{B}\in\mathbb{R}^{n\times n^2}$ and $\bm{C}\in\mathbb{R}^{n\times n^3}$ describe $2$nd- and $3$rd-order polynomial approximations of some nonlinear behavior. The matrix $\bm{D}\in\mathbb{R}^{n\times nm}$ captures the coupling between the state variables and input signals and $\bm{E}\in\mathbb{R}^{n\times m}$ describes how the input signals are injected into the circuit.
This differential equation will serve as the basis in the following model order reduction applications. 
\begin{figure}[t]
\centering
\includegraphics[width=3.1in]{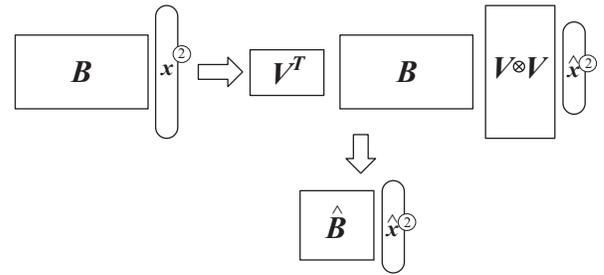}
\caption{Traditional projection-based nonlinear model order reduction methods reduce a large system matrix $\bm{B}$ to a small but dense matrix $\bm{\hat{B}}$ through an orthogonal projection matrix $\bm{V}$.} 
\label{fig:reducedB_flat}
\end{figure}

\begin{table*}[!t]
\centering
\caption{Computation and storage complexities of different nonlinear model order reduction approaches on a $q$-state reduced system with $d$th-order nonlinearity.}
\label{table:nmor} 
\begin{tabular}{lrrr}
Reduction methods & Function evaluation cost & Jacobian matrix evaluation cost &Storage cost \\ \midrule
Traditional matrix-based method~\cite{JR:99,LP:05,P:03,G:11,ZLWFW:12} & {$O(q^{d+1})$}  & {$O(q^{d+2})$} & {$O(q^{d+1})$}\\
Tensor-based method~\cite{LDW:15} & {$O(qdr)$}  & {$O(q^2dr)$} & {$O(qdr)$}\\
Symmetric tensor-based method~\cite{DLBKW:16} & {$O(qr)$}  & {$O(q^2r)$} &  {$O(qr)$}  
\end{tabular}
\end{table*}
\begin{figure*}[t]
\centering
\subfigure[]{\label{fig:tnode}
\includegraphics[width=2.5in]{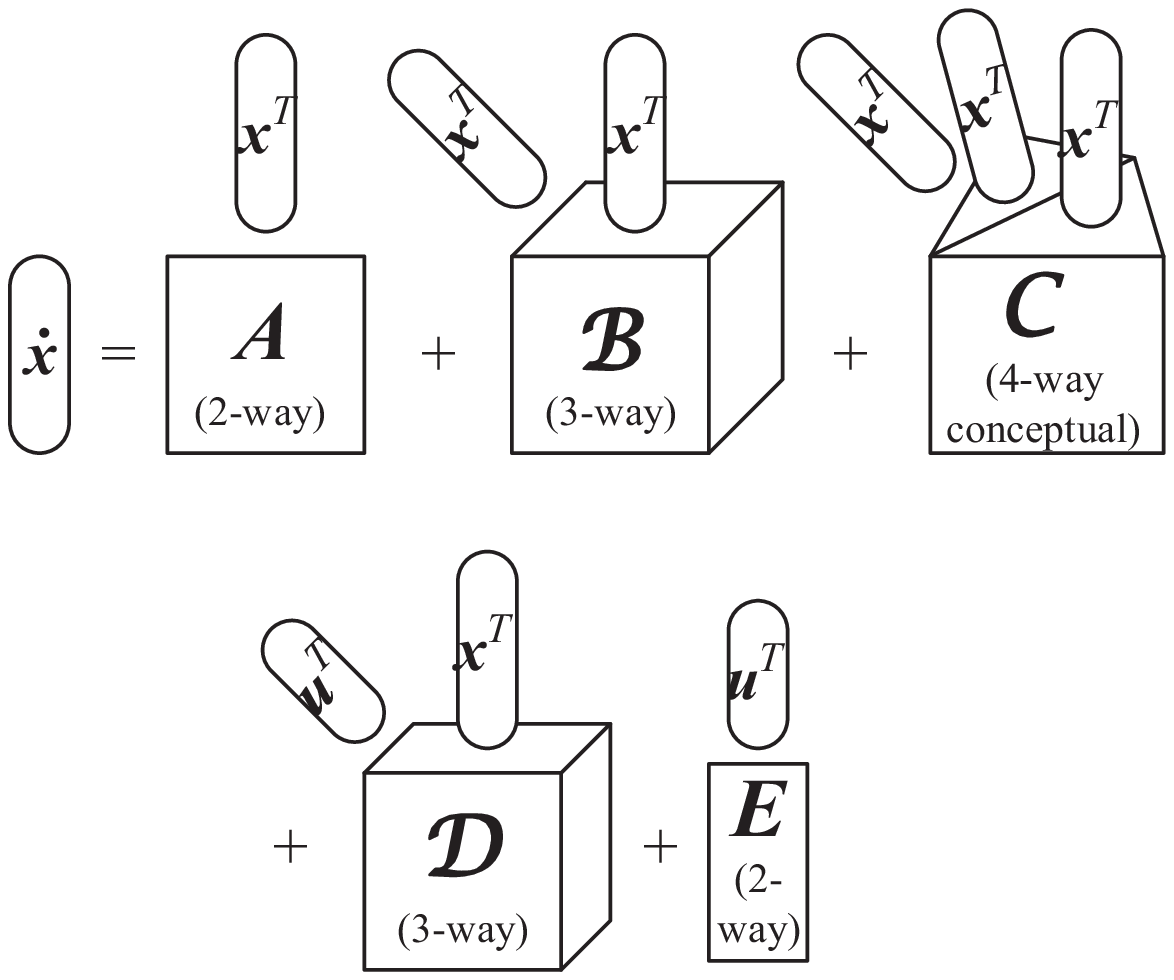}} \hspace{0.2in}
\subfigure[]{\label{fig:reducedb}
\includegraphics[width=3in]{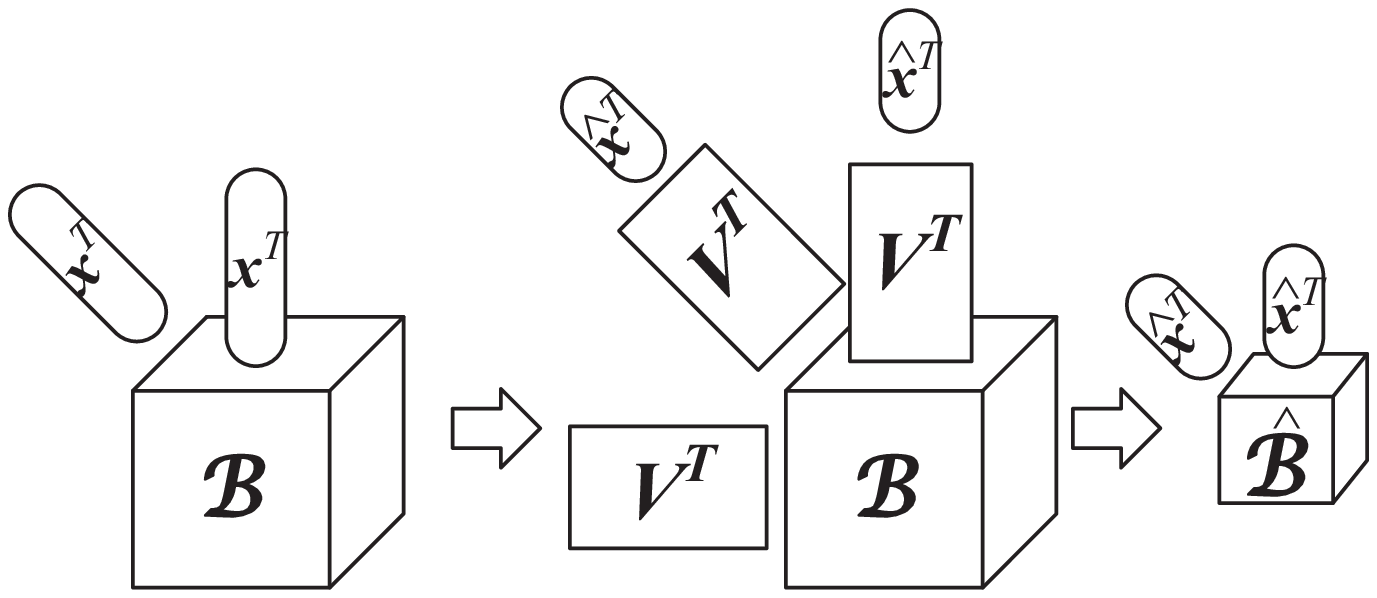}}
\caption{Tensor structures used in~\cite{LDW:15}. (a) Tensor representation of the original nonlinear system in~\eqref{eq:ode}; (b) tensor $\ten{B}$ is reduced to a compact tensor $\hat{\ten{B}}$ with a projection matrix $\bm{V}$ in~\cite{LDW:15}.}
\label{fig:tnmor}
\end{figure*}

\textbf{Matrix-based Nonlinear Model Order Reduction}. The idea of nonlinear model order reduction is to extract a compact reduced-order model that accurately approximates the input-output relationship of the original large nonlinear system. Simulation of the reduced-order model is usually much faster, so that efficient and reliable system-level verification is obtained. For instance, projection-based nonlinear model order reduction methods reduce the original system in~\eqref{eq:ode} to a compact reduced model with size $q\ll n$
\begin{align}
\bm{\dot{ \hat{x}}}=\bm{\hat{A}\hat{x}}+\bm{\hat{B}} \bm{\hat{x}}\kpr{2}+\bm{\hat{C}}\bm{\hat{x}}\kpr{3} + \bm{\hat{D}}(\bm{u} \otimes \bm{\hat{x}}) +\bm{\hat{E}u},
\label{eqn:romode}
\end{align}
where $\hat{x} \in \mathbb{R}^{q}$, $\bm{\hat{A}}\in\mathbb{R}^{q\times q}$, $\bm{\hat{B}}\in\mathbb{R}^{q\times q^2}$,  $\bm{\hat{C}}\in\mathbb{R}^{q\times q^3}$, $\bm{\hat{D}}\in\mathbb{R}^{q\times qm}$ and $\bm{\hat{E}}\in\mathbb{R}^{q\times m}$. The reduction is achieved through applying an orthogonal projection matrix $\bm{V}\in\mathbb{R}^{n\times q}$ on the system matrices in~\eqref{eq:ode}. Fig.~\ref{fig:reducedB_flat} illustrates how projection-based methods reduce $\bm{B}$ to a dense system matrix $\bm{\hat{B}}$ with a smaller size.

Most traditional matrix-based weakly nonlinear model order reduction methods~\cite{JR:99,LP:05,P:03,G:11,ZLWFW:12} suffer from the exponential growth of the size of the reduced system matrices $\bm{\hat{B}}, \bm{\hat{C}}, \bm{\hat{D}}$. As a result, simulating high-order strongly nonlinear reduced models is sometimes even slower than simulating the original system.

\textbf{Tensor-based Nonlinear Model Order Reduction}. A tensor-based reduction scheme was proposed in~\cite{LDW:15}. The coefficient matrices $\bm{B}, \bm{C}, \bm{D}$ of the polynomial system~\eqref{eq:ode} were reshaped into the respective tensors $\ten{B} \in \mathbb{R}^{n\times n\times n}$, $\ten{C} \in \mathbb{R}^{n\times n\times n \times n}$ and $\ten{D} \in \mathbb{R}^{n\times n\times m}$, as demonstrated in Fig.~\ref{fig:tnode}. These tensors were then decomposed via e.g. CPD, Tucker or Tensor Train rank-1 SVD, resulting in a tensor approximation of \eqref{eq:ode} as
\begin{align}\label{eqn:tenapprox}
\bm{\dot x} = & \bm{Ax} + [\![\bm{B}^{(1)}, \bm{x}^T\bm{B}^{(2)}, \bm{x}^T\bm{B}^{(3)}]\!]\nonumber\\
&+ [\![\bm{C}^{(1)}, \bm{x}^T\bm{C}^{(2)}, \bm{x}^T\bm{C}^{(3)}, \bm{x}^T\bm{C}^{(4)}]\!]\nonumber\\
&+[\![\bm{D}^{(1)}, \bm{x}^T\bm{D}^{(2)}, \bm{u}^T\bm{D}^{(3)}]\!]+\bm{Eu},
\end{align}
where $\bm{B}^{(k)},\bm{C}^{(k)},\bm{D}^{(k)},$ denote the $k$th-mode factor matrix from the polyadic decomposition of the tensors $\ten{B},\ten{C},\ten{D}$ respectively. Consequently, the reduced-order model inherits the same tensor structure as~\eqref{eqn:tenapprox} (with smaller sizes of the mode factors). If we take tensor $\ten{B}$ as an example, its reduction process in~\cite{LDW:15} is shown in Fig.~\ref{fig:reducedb}.

\textbf{Computational and Storage Benefits}. Unlike previous matrix-based approaches, simulation of the tensor-structure reduced model completely avoids the overhead of solving high-order dense system matrices, since the dense Kronecker products in~\eqref{eq:ode} are resolved by matrix-vector multiplications between the mode factor matrices and the state vectors. Therefore, substantial improvement on efficiency can be achieved. Meanwhile, these mode factor matrices can significantly reduce the memory requirement since they replace all dense tensors and can be reduced and stored beforehand. Table~\ref{table:nmor} shows the computational complexities of function and Jacobian matrix evaluations when simulating a reduced model with $d$th-order nonlinearity, where $r$ denotes the tensor rank used in the polyadic decompositions in~\cite{LDW:15}. The storage costs of those methods are also listed in the last column of Table~\ref{table:nmor}.


\textbf{Symmetric Tensor-based Nonlinear Model Order Reduction}. A symmetric tensor-based order reduction method in~\cite{DLBKW:16} further utilizes the all-but-first-mode partial symmetry of the system tensor $\ten{B}$ ($\ten{C}$), i.e., the mode factors of $\ten{B}$ ($\ten{C}$) are exactly the same, except for the first mode only. This partial symmetry property is also kept by its reduced-order model. The symmetric tensor-based reduction method in~\cite{DLBKW:16} provides further improvements of computation performance and storage requirement over~\cite{LDW:15}, as shown in the last row of Table~\ref{table:nmor}.

\subsection{Volterra-Based Simulation and Identification for ICs}

Volterra theory has long been used in analyzing communication systems and in nonlinear control~\cite{Bed:71,R:81}. It can be regarded as a kind of Taylor series with memory effects since its evaluation at a particular time point requires input information from the past. Given a certain input and a black-box model of a nonlinear system with time/frequency-domain Volterra kernels, the output response can be computed by the summation of a series of multidimensional convolutions. For instance, a $3$rd-order response can be written in a discretized form as
\begin{align}\label{eqn:conv3rd}
y_3[k]=\sum^M_{m_1=1}\sum^M_{m_2=1}\sum^M_{m_3=1}h_3[m_1,m_2,m_3]\prod^3_{i=1} u[k-m_i],
\end{align}
where $h_3$ denotes the $3$rd-order Volterra kernel, $u$ is the discretized input and $M$ is the memory. Such a multidimensional convolution is usually done by multidimensional fast Fourier transforms (FFT) and inverse fast Fourier transforms (IFFT). Although the formulation does not preclude itself from modeling strong nonlinearities, the exponential complexity growth in multidimensional FFT/IFFT computations results in the curse of dimensionality that forbids its practical implementation.

\begin{figure}[!t]
	\centering
		\includegraphics[width=2.8in]{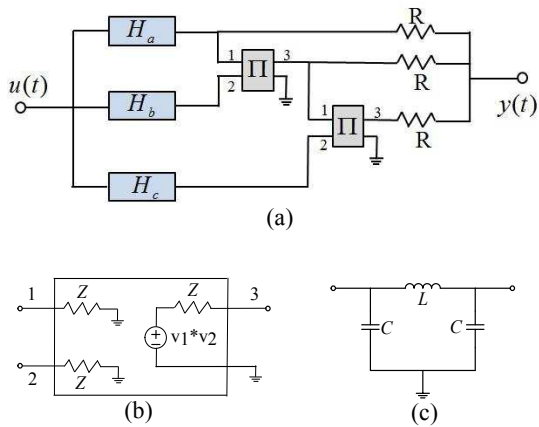} 
\caption{(a) System diagram of a $3$rd-order mixer circuit. The symbol $\Pi$ denotes a mixer; (b) the equivalent circuit of the mixer. $Z=R=50\,\Omega$; (c) the circuit schematic diagram of the low-pass filters $H_a$, $H_b$ and $H_c$, with $L=42.52\,\textrm{nH}$ and $C=8.5\,\textrm{pF}.$}
	\label{fig:mixer}
\end{figure}

\textbf{Tensor-Volterra Model-based Simulation}. Obviously, the $3$rd-order Volterra kernel $h_3$ itself can be viewed as a $3$-way tensor. 
By compressing the Volterra kernel into a polyadic decomposition, it is proven in~\cite{LXBJDW:15} that the computationally expensive multidimensional FFT/IFFT can be replaced by a number of cheap one-dimensional FFT/IFFTs without compromising much accuracy. 

\textbf{Computational and Storage Benefits}. The chosen rank for the polyadic decomposition has a significant impact on both the accuracy and efficiency of the simulation algorithm. In~\cite{LXBJDW:15}, the ranks were chosen a priori and it was demonstrated that the computational complexity for the tensor-Volterra based method to calculate an $d$th-order response is in $O((R_{\text{real}}+R_{\text{imag}})dm\log m)$, where $m$ is the number of steps in the time/frequency axis, and $R_{\textrm{real}}$ and $R_{\textrm{imag}}$ denote the prescribed ranks of the polyadic decomposition used for the real and imaginary parts of the Volterra kernel, respectively. In contrast, the complexity for the traditional multidimensional FFT/IFFT approach is in $O(dm^d\log m)$. In addition, the tensor-Volterra model requires the storage of only the factor matrices in memory, with space complexity $O((R_{\text{real}}+R_{\text{imag}})dm)$, while $O(m^d)$ memory is required for the conventional approach.

\begin{figure*}[t]
\centering
\subfigure[]{\label{fig:mixer_t}
\includegraphics[width=2in, height=2.2in]{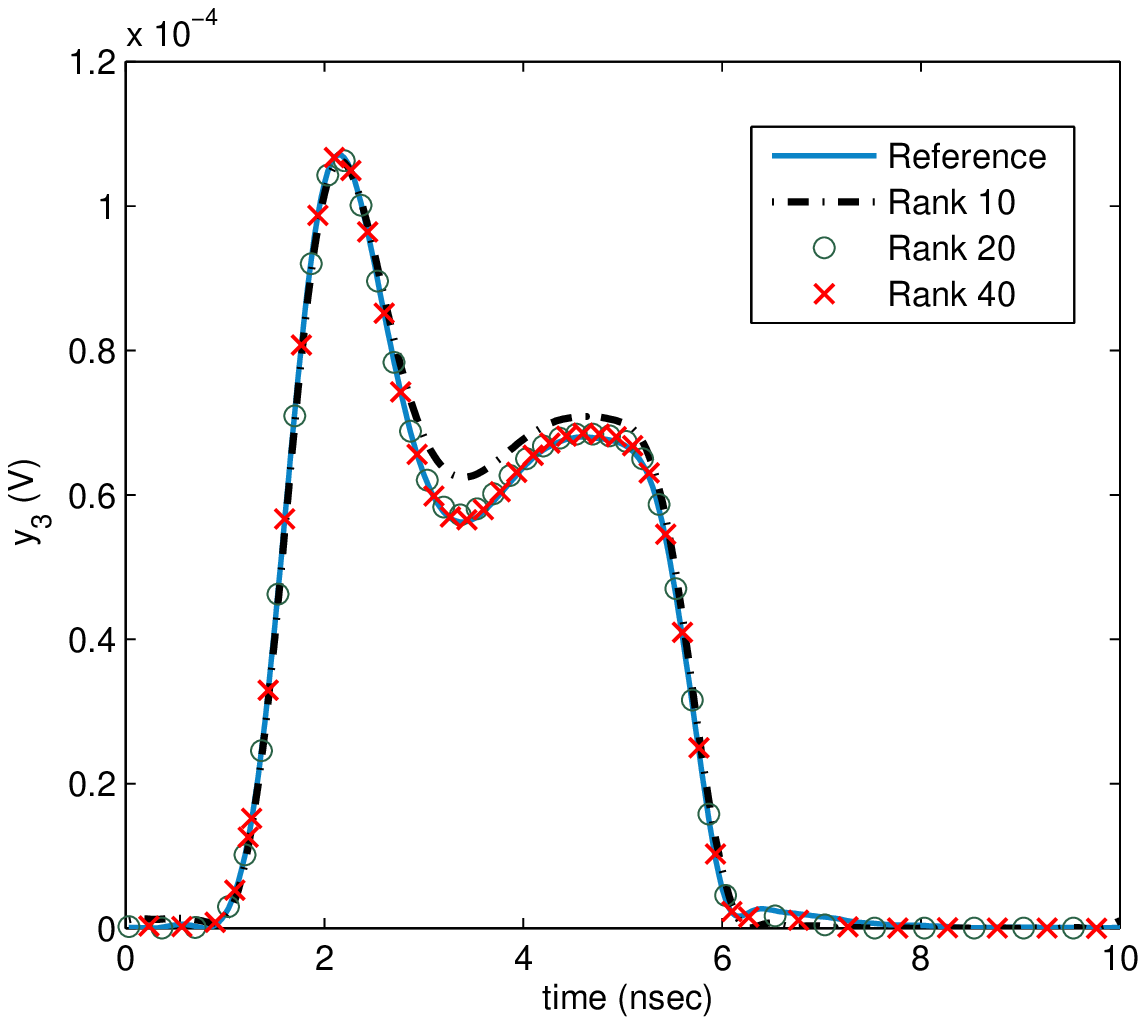}} \hspace{0.2in}
\subfigure[]{\label{fig:mixer_terr}
\includegraphics[width=2in, height=2.2in]{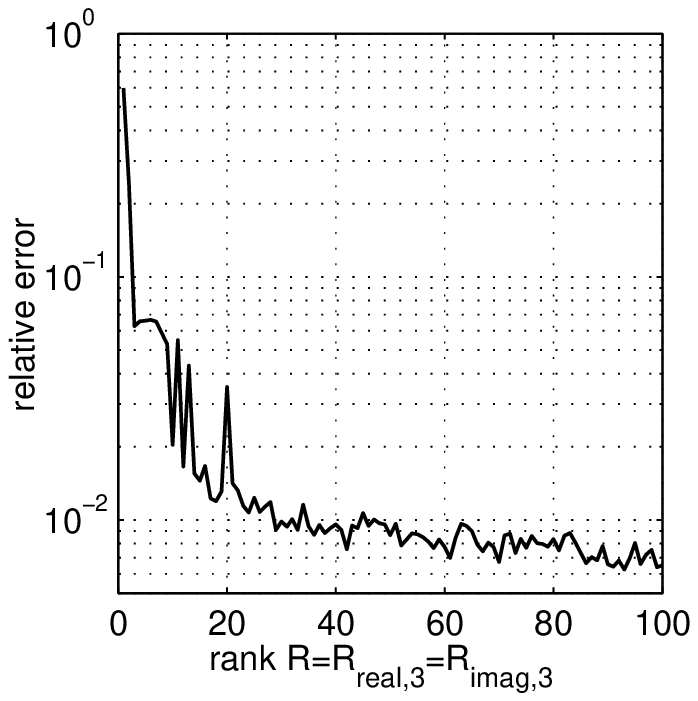}} \hspace{0.2in}
\subfigure[]{\label{fig:mixer_tspeed}
\includegraphics[width=2in, height=2.2in]{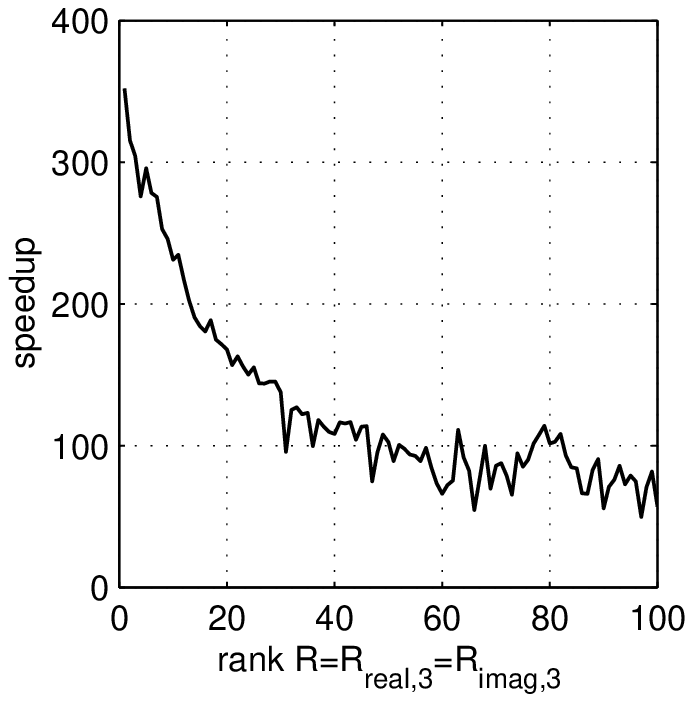}}
\caption{Numerical results of the mixer. (a) Time-domain results of $y_3$ computed by the method in~\cite{LXBJDW:15} with different rank approximations; (b) relative errors of~\cite{LXBJDW:15} with different ranks; (c) speedups brought by~\cite{LXBJDW:15} with different ranks.}
\label{fig:mixer_result}
\end{figure*}

In~\cite{LXBJDW:15}, the method was applied to compute the time-domain response of a $3$rd-order mixer system shown in Fig.~\ref{fig:mixer}. The $3$rd-order response $y_3$ is simulated to a square pulse input with $m=201$ time steps. As shown in Fig.~\ref{fig:mixer_result}(a), a rank-20 (or above) polyadic decomposition for both the real and imaginary parts of the kernel tensor matched the reference result from multidimensional FFT/IFFT fairly well. Figs.~\ref{fig:mixer_result}(b) and (c) demonstrate a certain trade-off between the accuracy and efficiency when using different ranks for the polyadic decomposition. Nonetheless, a 60x speedup is still achievable for ranks around 100 with a 0.6\% error.

\textbf{System Identification}. In~\cite{FKB:12,AF:04,iccad16volterra}, similar tensor-Volterra models were used to identify the black-box Volterra kernels $h_i$. It was reported in~\cite{FKB:12,AF:04,iccad16volterra} that given certain input and output data, identification of the kernels in the polyadic decomposition form could significantly reduce the parametric complexity with good accuracy.

\section{Future Topics: EDA Applications}
\label{sec:future_app}
This section describes some EDA problems that could be potentially solved with, or that could benefit significantly from employing tensors. Since many EDA problems are characterized by high dimensionality, the potential application of tensors in EDA can be vast and is definitely not limited to the topics summarized below.

\subsection{EDA Optimization with Tensors}
\label{subsec:opt}
Many EDA problems require solving a large-scale optimization problem in the following form:
 \begin{align}
\label{EDA_opt}
& \min_{\bm{x}} \; f(\bm{x}), \;\;\; {\rm s.}\; {\rm t.} \;\;\;\bm{x} \in {\cal C}
\end{align}
where $\bm{x}=[x_1, \cdots, x_n]$ denotes $n$ design or decision variables, $f(\bm{x})$ is a cost function (e.g., power consumption of a chip, layout area, signal delay), and $ {\cal C}$ is a feasible set specifying some design constraints. This formulation can describe problems such as circuit optimization~\cite{hershenson2001optimal,li2004robust, xu2005opera}, placement~\cite{Shahookar1991}, routing~\cite{cong1996performance}, and power management~\cite{benini1999policy}. The optimization problem (\ref{EDA_opt}) is computationally expensive if $\bm{x}$ has many elements.

It is possible to accelerate the above large-scale optimization problems by exploiting tensors. By adding some extra variables $\hat{\bm{x}}$ with $\hat {n}$ elements, one could form a longer vector $\bar{\bm{x}}=[\bm{x}, \hat{\bm{x}}]$ such that $\bar{\bm{x}}$ has $n_1\times \cdots \times n_d$ variables in total.  Let $\bar{\ten{X}}$ be a tensor such that $\bar{\bm{x}}=\textrm{vec}(\bar{\ten{X}})$, let $\bm{x}=\bm{Q}\bar{\bm{x}}$ with $\bm{Q}$ being the first $n$ rows of an identity matrix, then (\ref{EDA_opt}) can be written in the following tensor format:
 \begin{align}
\label{EDA_opt_ten}
& \min_{\bar{\ten{X}}} \; \bar{f}(\bar{\ten{X}}), \;\;\; {\rm s.}\; {\rm t.} \;\;\;\bar{\ten{X}} \in \bar{{\cal C}}
\end{align}
with $ \bar{f}(\bar{\ten{X}})=f(\bm{Q} \textrm{vec}(\bar{\ten{X}}))$ and $\bar{{\cal C}}=\{ \bar{\ten{X}}| \bm{Q} \textrm{vec}(\bar{\ten{X}}) \in {\cal C} \}$.

Although problem~\eqref{EDA_opt_ten} has more unknown variables than~\eqref{EDA_opt}, the low-rank representation of tensor~$\bar{\ten{X}}$ may have much fewer unknown elements. Therefore, it is highly possible that solving \eqref{EDA_opt_ten} will require much lower computational cost for lots of applications. 

\subsection{High-Dimensional Modeling and Simulation}

Consider a general algebraic equation resulting from a high-dimensional modeling or simulation problem
\begin{equation}
\label{EDA_lieq}
\bm{g}(\bm{x})=0, \;\;\;{\rm with} \; \bm{x} \in \mathbb{R}^N \;{\rm and} \;  N=n_1 \times n_2 \cdots \times n_d
\end{equation}
which can be solved by Newton's iteration. For simplicity, we assume $n_i=n$. When an iterative linear equation solver is applied inside a Newton's iteration, it is possible to solve this problem at the complexity of $O(N)=O(n^d)$. However, since $N$ is an exponential function of $n$, the iterative matrix solver quickly becomes inefficient as $d$ increases. Instead, we rewrite \eqref{EDA_lieq} as the following equivalent optimization problem:
\begin{equation}
\label{EDA_lieq_opt}
\min_{\bm{x}} \; f(\bm{x})=\| \bm{g}(\bm{x})\|_2^2, \;\;\; {\rm s.}\; {\rm t.} \;\;\;\bm{x} \in \mathbb{R}^N. \nonumber
\end{equation}
This least-square optimization is a special case of (\ref{EDA_opt}), and thus the tensor-based optimization idea may be exploited to solve the above problem at the cost of $O(n)$.

A potential application lies in the PDE or integral equation solvers for device simulation. Examples include the Maxwell equations for parasitic extraction~\cite{FastCap,FastHenrry,pFFT}, the Navier-Stokes equation describing bio-MEMS~\cite{VasilyevRW06}, and the Poisson equation describing heating effects~\cite{YuZYQ13}. These problems can be described as (\ref{EDA_lieq}) after numerical discretization. The tensor representation of $\bm{x}$ can be easily obtained based on the numerical discretization scheme. For instance, on a regular 3-D cubic structure, a finite-difference or finite-element discretization may use $n_x, n_y$ and $n_z$ discretization elements in the $x$, $y$ and $z$ directions respectively. Consequently, $\bm{x}$ could be compactly represented as a $3$-way tensor with size $n_x \times n_y \times n_z$ to exploit its low-rank property in the spatial domain. 

This idea can also be exploited to simulate multi-rate circuits or multi-tone RF circuits. In both cases, the tensor representation of $\bm{x}$ can be naturally obtained based on the time-domain discretization or multi-dimensional Fourier transform. In multi-tone harmonic balance~\cite{Melville:mt1995, Carvalho:mt1998}, the dimensionality $d$ is the total number of RF inputs. In the multi-time PDE solver~\cite{Jaijeet:PDE2001}, $d$ is the  number of time axes describing different time scales.

\subsection{Process Variation Modeling}
\begin{figure}[t]
	\centering
		\includegraphics[width=2.8in]{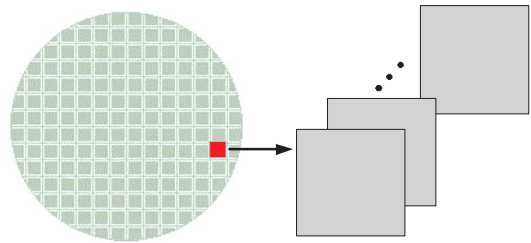} 
\caption{Represent multiple testing chips on a wafer as a single tensor. Each slice of the tensor captures the spatial variations on a single die. }
	\label{fig:die_wafer}
\end{figure}

In order to characterize the inter-die and intra-die process variations across a silicon wafer with $k$ dice, one may need to measure each die with an $m\times n$ array of devices or circuits~\cite{variation2008, yu2014remembrance, zhang2011virtual}. The variations of a certain parameter (e.g., transistor threshold voltage) on the $i$th die can be described by matrix $\bm{A}_i \in \mathbb{R}^{m\times n}$, and thus one could describe the whole-wafer variation by stacking all matrices into a tensor $\ten{A} \in \mathbb{R}^{k\times m \times n}$, with $\bm{A}_i$ being the $i$th slice. This representation is graphically shown in Fig.~\ref{fig:die_wafer}.

Instead of measuring each device on each die (which requires $kmn$ measurements in total), one could measure only a few devices on each wafer, then estimate the full-wafer variations using tensor completion. One may employ convex optimization to locate the globally optimal solution of this $3$-way tensor completion problem. 

\section{Future Topics: Theoretical Challenges}
\label{sec:future}
Tensor theory is by itself an active research topic. This section summarizes some theoretical open problems. 
\subsection{Challenges in Tensor Decomposition} 
Polyadic and tensor train decompositions are preferred for high-order tensors due to their better scalability. In spite of their better computational scalability, the following challenges still exist:
\begin{itemize}
\item \textbf{Rank Determination in CPD.} The tensor ranks are usually determined by two methods. First, one may fix the rank and search for the tensor factors. Second, one may increase the rank incrementally to achieve an acceptable accuracy. Neither methods are optimal in the theoretical sense.

\item \textbf{Optimization in Polyadic Decomposition.} Most rank-$r$ polyadic decomposition algorithms employ alternating least-squares (ALS) to solve non-convex optimization problems. Such schemes do not guarantee the global optimum, and thus it is highly desirable to develop global optimization algorithms for the CPD.

\item \textbf{Faster Tensor Train Decomposition.} Computing the tensor train decomposition requires the computation of many low-rank decompositions. The state-of-the-art implementation employs ``cross approximation'' to perform low-rank approximations~\cite{ttcross}, but it still needs too many iterations to find a ``good'' representation.

\item \textbf{Preserving Tensor Structures and/or Properties.} In some cases, the given tensor may have some special properties such as symmetry or non-negativeness. These properties need to be preserved in their decomposed forms for specific applications.
\end{itemize}

\subsection{Challenges in Tensor Completion} 
Major challenges of tensor completion include:
\begin{itemize}
\item \textbf{Automatic Rank Determination.} In high-dimensional tensor completion, it is important to determine the tensor rank automatically. Although some probabilistic approaches such as variational Bayesian methods~\cite{Zhao:2015, Zhao:arxiv2015} have been reported, they are generally not robust for very high-order tensors.

\item \textbf{Convex Tensor Completion.} Most tensor completion problems are formulated as non-convex optimization problems. Nuclear-norm minimization is convex, but it is only applicable to low-order tensors. Developing a scalable convex formulation for the minimal-rank completion still remains an open problem for high-order cases.

\item \textbf{Robust Tensor Completion.} In practical tensor completion, the available tensor elements from measurement or simulations can be noisy or even wrong. For these problems, the developed tensor completion algorithms should be robust against noisy input. 

\item \textbf{Optimal Selection of Samples.} Two critical fundamental questions should be addressed. First, how many samples are required to (faithfully) recover a tensor? Second, how can we select the samples optimally?
\end{itemize}

\section{Conclusion}
\label{sec:conclusion}
By exploiting low-rank and other properties of tensors (e.g., sparsity, symmetry), the storage and computational cost of many challenging EDA problems can be significantly reduced. For instance, in the high-dimensional stochastic collocation modeling of a CMOS ring oscillator, exploiting tensor completion required only a few hundred circuit/device simulation samples vs. the huge number of simulations (e.g., $10^{27}$) required by standard approaches to build a stochastic model of similar accuracy. When applied to hierarchical uncertainty quantification, a tensor-train approach allowed the easy handling of an extremely challenging MEMS/IC co-design problem with over 180 uncorrelated random parameters describing process variations. In nonlinear model order reduction, the high-order nonlinear terms were easily approximated by a tensor-based projection framework. Finally, a $60\times$ speedup was observed when using tensor computation in a 3rd-order Volterra-series nonlinear modeling example, while maintaining a $0.6\%$ relative error compared with the conventional FFT/IFFT approach. These are just few initial representative examples for the huge potential that a tensor computation framework can offered to EDA algorithms. We believe that the space of EDA applications that could benefit from the use of tensors is vast and remains mainly unexplored, ranging from EDA optimization problems, to device field solvers, and to process variation modeling.

\section*{Acknowledgement}
This work was supported in part by the NSF NEEDS program, the AIM Photonics Program, the Hong Kong Research Grants Council under Project 17212315, the University Research Committee of The University of Hong Kong, and the MIT Greater China Fund for Innovation.

\appendices 

\section{Additional Notations and Definitions}
\label{appd:notations}
\textbf{Diagonal, Cubic and Symmetric Tensors}. The diagonal entries of a tensor $\ten{A}$ are the entries $a_{i_1i_2\cdots i_d}$ for which $i_1=i_2=\cdots = i_d$. A tensor $\ten{S}$ is diagonal if all of its non-diagonal entries are zero. A cubical tensor is a tensor for which $n_1=n_2=\cdots=n_d$. A cubical tensor $\ten{A}$ is symmetric if $a_{i_1\cdots i_d}=a_{\pi(i_1,\ldots,i_d)}$ where $\pi(i_1,\ldots,i_d)$ is any permutation of the indices.

The \textbf{Kronecker product}~\cite{VL:00} is denoted by $\otimes$. We use the notation $x\kpr{d}= x\otimes x\otimes \cdots \otimes x$ for the $d$-times repeated Kronecker product.

\definition{\textbf{Reshaping.}} Reshaping, also called \textbf{unfolding}, is another often used tensor operation. The most common reshaping is the matricization, which reorders the entries of $\ten{A}$ into a matrix. The mode-$n$ matricization of a tensor $\ten{A}$, denoted $\ten{A}_{(n)}$, rearranges the entries of $\ten{A}$ such that the rows of the resulting matrix are indexed by the $n$th tensor index $i_n$. The remaining indices are grouped in ascending order.

\example The $3$-way tensor of Fig.~\ref{fig:tensorexample} can be reshaped as a $2\times 12$ matrix or a $3\times 8$ matrix, and so forth. The mode-$1$ and mode-$3$ unfoldings are
\begin{gather*}
\ten{A}_{(1)} =
\begin{pmatrix}
1& 4& 7& 10& 13& 16& 19& 22\\
2& 5& 8& 11& 14& 17& 20& 23\\
3& 6& 9& 12& 15& 18& 21& 24\\
\end{pmatrix},\\
\ten{A}_{(3)} =
\begin{pmatrix}
1 & 2 & 3 & 4 & \cdots & 9 & 10& 11 & 12\\
13& 14& 15& 16 & \cdots & 21 & 22 & 23 & 24
\end{pmatrix}.
\end{gather*}
The column indices of $\ten{A}_{(1)},\ten{A}_{(3)}$ are $[i_2i_3]$ and $[i_1i_2]$, respectively. 

\definition{\textbf{Vectorization.}} Another important reshaping is the vectorization. The vectorization of a tensor $\ten{A}$, denoted $\textrm{vec}(\ten{A})$, rearranges its entries in one vector. 

\example For the tensor in Fig.~\ref{fig:tensorexample}, we have 
\begin{align*}
\textrm{vec}(\ten{A}) \;=\; \begin{pmatrix}
1& 2& \cdots &24
\end{pmatrix}^T.
\end{align*}

\section{Computation and Variants of the Polyadic Decomposition}
\label{app:cpvariant}
\textbf{Computing Polyadic Decompositions}. Since the tensor rank is not known a priori, in practice, one usually computes a low-rank $r < R$ approximation of a given tensor $\ten{A}$ by minimizing the Frobenius norm of the difference between $\ten{A}$ and its approximation. Specifically, the user specifies $r$ and then solves the minimization problem
$$\argmin{\ten{D},\bm{U}^{(1)},\ldots,\bm{U}^{(d)}}||\ten{A}-[\![\ten{D};\bm{U}^{(1)}, \ldots, \bm{U}^{(d)}]\!]||_F
$$
where $\ten{D} \in \mathbb{R}^{r \times r \times \cdots \times r}, \bm{U}^{(i)} \in \mathbb{R}^{n_i \times r} (i=\{1,\ldots,d\})$. One can then increment $r$ and compute new approximations until a ``good enough'' fit is obtained. A common method for solving this optimization problem is the Alternating Least Squares (ALS) method~\cite{harshman1970fpp}. Other popular optimization algorithms are nonlinear conjugate gradient methods, quasi-Newton or nonlinear least squares (e.g. Levenberg-Marquardt)~\cite{Sorber2013}. The computational complexity per iteration of the ALS, Levenberg-Marquardt (LM) and Enhanced Line Search (ELS) methods to compute a polyadic decomposition of a $3$-way tensor, where $n=\textrm{min}(n_1,n_2,n_3)$, are given in Table~\ref{table:tdecompcomp}.
\begin{table}[ht]
\centering
\caption{Computational costs of 3 tensor decomposition methods for a $3$-way tensor~\cite{Comon2009}.}
\label{table:tdecompcomp} 
\begin{tabular}{@{}lr@{}}
Methods & Cost per iteration \\ \midrule
\\
ALS & {$(n_2n_3+n_1n_3+n_1n_2)(7n^2+n)+3nn_1n_2n_3$}   \\
\\
LM & {$n_1n_2n_3(n_1+n_2+n_3)^2n^2$}   \\
\\
ELS & {$(8n+9)n_1n_2n_3$}  
\end{tabular}
\end{table}

Two variants of the polyadic decomposition are summarized below.
\subsubsection{PARATREE or tensor-train rank-1 SVD (TTr1SVD)} This polyadic decomposition~\cite{ttr1svd,salmi2009sequential} consists of orthogonal rank-1 terms and is computed by consecutive reshapings and SVDs. This computation implies that the obtained decomposition does not need an initial guess and will be unique for a fixed order of indices. Similar to SVD in the matrix case, this decomposition has an approximation error easily expressed in terms of the $\sigma_i$'s~\cite{ttr1svd}. 

\subsubsection{CPD for Symmetric Tensors} The CPD of a symmetric tensor does not in general result in a summation of symmetric rank-1 terms. In some applications, it is more meaningful to enforce the symmetric constraints explicitly, and write $\ten{A}=\sum_{i=1}^R \lambda_i \bm{v}_i^d$, where $\lambda_i \in \mathbb{R}, \ten{A}$ is a $d$-way symmetric tensor. Here $\bm{v}_i^d$ is a shorthand for the $d$-way outer product of a vector $\bm{v}_i$ with itself, i.e., $\bm{v}^d_i=\bm{v}_i \circ \bm{v}_i \circ \cdots \circ \bm{v}_i$.

\section{Higher-order SVD}
\label{app:hosvd}
The Higher-Order SVD (HOSVD)~\cite{Lathauwer:HOSVD} is obtained from the Tucker decomposition when the factor matrices $\bm{U}^{(i)}$ are orthogonal, when any two slices of the core tensor $\ten{S}$ in the same mode are orthogonal, $\langle \ten{S}_{i_k=p},\ten{S}_{i_k=q} \rangle = 0$ if $p \neq q$ for any $k=1,\ldots,d$, and when the slices of the core tensor $\ten{S}$ in the same mode are ordered according to their Frobenius norm, $||\ten{S}_{i_k=1}|| \geq ||\ten{S}_{i_k=2}|| \geq \cdots \geq ||\ten{S}_{i_k=n_k}||$ for $k=\{1,\ldots,d\}$.
Its computation consists of $d$ SVDs to compute the factor matrices and a contraction of their inverses with the original tensor to compute the HOSVD core tensor. For a 3-way tensor this entails a computational cost of $2n_1n_2n_3(n_1 +n_2 +n_3)+5(n_1^2 n_2n_3 +n_1n_2^2n_3 +n_1n_2n_3^2 )â2(n_1^3 +n_2^3 +n_3^3 )/3â(n_1^3 +n_2^3 +n_3^3 )/3$~\cite{Comon2009}.

\bibliographystyle{IEEEtran}
\bibliography{IEEEabrv,keynote}
\begin{IEEEbiography}[]{Zheng Zhang} (M'15) received the Ph.D degree (2015) in Electrical Engineering and Computer Science from the Massachusetts Institute of Technology (MIT), Cambridge, MA. Currently he is a Postdoc Associate with the Research Laboratory of Electronics at MIT. His research interests include uncertainty quantification, tensor and model order reduction, with diverse engineering applications including nanoelectronics, energy systems and biomedical imaging. His industrial experiences include Coventor Inc. and Maxim-IC; academic visiting experiences include UC San Diego, Brown University and Politechnico di Milano; government lab experiences include the Argonne National Laboratory.

Dr. Zhang received the 2016 ACM Outstanding Ph.D Dissertation Award in Electronic Design Automation, the 2015 Doctoral Dissertation Seminar Award (i.e., Best Thesis Award) from the Microsystems Technology Laboratory of MIT, the 2014 Donald O. Pederson Best Paper Award from IEEE Transactions on Computer-Aided Design of Integrated Circuits and Systems, the 2014 Chinese Government Award for Outstanding Students Abroad, and the 2011 Li Ka-Shing Prize from the University of Hong Kong. 
\end{IEEEbiography}

\begin{IEEEbiography}[]{Kim Batselier} (M'13) received the M.S. degree in Electro-Mechanical Engineering and the Ph.D. Degree in Engineering Science from the KULeuven, Belgium, in 2005 and 2013 respectively. He worked as a research engineer at BioRICS on automated performance monitoring until 2009. He is currently a Post-Doctoral Research Fellow at The University of Hong Kong since 2013. His current research interests include linear and nonlinear system theory/identification, algebraic geometry, tensors, and numerical algorithms.
\end{IEEEbiography}

\begin{IEEEbiography}[]{Haotian Liu} (S'11) received the B.S. degree in Microelectronic Engineering from Tsinghua University, Beijing, China, in 2010, and the Ph.D. degree in Electronic Engineering from the University of Hong Kong, Hong Kong, in 2014. He is currently a software engineer with Cadence Design Systems, Inc. San Jose, CA.

In 2014, Dr. Liu was a visiting scholar with the Massachusetts Institute of Technology (MIT), Cambridge, MA. His research interests include numerical simulation methods for very-large-scale integration (VLSI) circuits, model order reduction, parallel computation and static timing analysis.
\end{IEEEbiography}

\begin{IEEEbiography}[]{Luca Daniel} (S'98-M'03) is a Full Professor in the Electrical Engineering and Computer Science Department of the Massachusetts Institute of Technology (MIT). He received the Ph.D. degree in Electrical Engineering from the University of California, Berkeley, in 2003. Industry experiences include HP Research Labs, Palo Alto (1998) and Cadence Berkeley Labs (2001).

Dr. Daniel current research interests include integral equation solvers, uncertainty quantification and parameterized model order reduction, applied to RF circuits, silicon photonics, MEMs, Magnetic Resonance Imaging scanners, and the human cardiovascular system.

Prof. Daniel was the recipient of the 1999 IEEE Trans. on Power Electronics best paper award; the 2003 best PhD thesis awards from the Electrical Engineering and the Applied Math departments at UC Berkeley; the 2003 ACM Outstanding Ph.D. Dissertation Award in Electronic Design Automation; the 2009 IBM Corporation Faculty Award; the 2010 IEEE Early Career Award in Electronic Design Automation; the 2014 IEEE Trans. On Computer Aided Design best paper award; and seven best paper awards in conferences.
\end{IEEEbiography}

\begin{IEEEbiography}[]{Ngai Wong} (S'98-M'02) received his B.Eng. and Ph.D. degrees in Electrical and Electronic Engineering from The University of Hong Kong, Hong Kong, in 1999 and 2003, respectively.

Dr. Wong was a visiting scholar with Purdue University, West Lafayette, IN, in 2003. He is currently an Associate Professor with the Department of Electrical and Electronic Engineering, The University of Hong Kong. His current research interests include linear and nonlinear circuit modeling and simulation, model order reduction, passivity test and enforcement, and tensor-based numerical algorithms in electronic design automation (EDA).
\end{IEEEbiography}
\end{document}